\documentclass{amsproc}
\usepackage{amsmath,amssymb,latexsym}
\usepackage{amscd}
\usepackage{rotating}
\usepackage[mathscr]{eucal}

\def\ca{\operatorname{ca}}

\def\w{\tilde}
\def\h{\hat}

\def\amp{\operatorname{amp}}

\def\lim{\operatorname{lim}}
\def\Lef{\operatorname{Lef}}
\def\Leff{\operatorname{Leff}}

\def\codim{\operatorname{codim}}

\def\id{\operatorname{id}}

\def\Im{\operatorname{Im}}

\newtheorem{theorem}{Theorem}[section]

\newtheorem{corollary}[theorem]{Corollary}

\newtheorem{prop}[theorem]{Proposition}

\newtheorem{definition}[theorem]{Definition}
\newtheorem{rem}[theorem]{Remark}
\newtheorem{rems}[theorem]{Remarks}

\theoremstyle{definition}

\theoremstyle{remark}

\numberwithin{equation}{section}

\begin{document}

\title{On a Connectedness Theorem of Debarre}

\author{Lucian B\u adescu}
\address{Dipartmento di Matematica, Universit\`a degli Studi di Genova, Via Dodecaneso
35, 16146 Genova, Italy}
\email{badescu@dima.unige.it}


\date{October 15}

\dedicatory{To Andrew Sommese on the occasion of his sixtieth anniversary.}


\begin{abstract}
Under a slightly stronger hypothesis, one improves a connectedness result of Debarre \cite{De} for a product of two projective spaces in terms of the extension problem of formal-rational functions (see Theorems \ref{main} and \ref{main1} of the Introduction).\end{abstract}

\maketitle



\section{Introduction}

We start by recalling a basic definition regarding the extension of formal-rational functions. Given an irreducible algebraic variety $X$ over $k$ and a closed subvariety $Y$ of $X$, we
shall denote by $K(X)$ the field of rational functions of $X$, by $X_{/Y}$ the formal completion of
$X$ along $Y$, and by $K(X_{/Y})$ the $k$-algebra of formal-rational functions of $X$ along $Y$,
i.e. the global sections of the sheaf of total fractions of the structural sheaf $\mathscr O_{X_{/Y}}$
(see \cite{HM}, or \cite{Ha}, or also \cite{B1}, Chapter 9). There is a natural (injective) homomorphism of 
$k$-algebras $\alpha_{X,Y}:K(X)\to K(X_{/Y})$.

\begin{definition}\label{oo}  {\em Let $X$ be a complete irreducible variety over an algebraically closed
field $k$, and let $Y$ be a closed subvariety of $X$. 
According to Hironaka and Matsumura \cite{HM} we say that $Y$ is $G3$ in $X$ if the canonical map
$\alpha_{X,Y}:K(X)\to K(X_{/Y})$ is an isomorphism of $k$-algebras. In other words, $Y$ is $G3$ in $X$ if every formal rational-function of $X$ along $Y$ extends to a rational function on $X$.}\end{definition}

Before stating our main results we fix some notation that will be used throughout this paper.  Let $n_1,\ldots,n_s$ be $s$ positive integers ($s\geq 1$) and let $P:=\mathbb P^{n_1}\times\cdots\times\mathbb P^{n_s}$ the product of the projective spaces $\mathbb P^{n_1},\ldots,\mathbb P^{n_s}$ over 
$k$. For every non-empty subset $J$ of $I:=\{1,\ldots,s\}$ set $P_J:=\prod_{i\in J}\mathbb P^{n_i}$, and denote by $p_J$ the canonical projection of $P=P_I$ onto $P_J$ (so that $p_I\colon P\to P$ is the identity) and  by $\Delta$ the diagonal of $P\times P$. 

The starting point of this paper is the following connectedness result of Debarre \cite{De} (which generalizes the well-known connectedness theorem of Fulton--Hansen \cite{FH} to the case of a product of projective spaces):

\begin{theorem}[Debarre \cite{De}]\label{debarre} Under the above notation, let $f\colon X\to P\times P$ be a morphism from a complete irreducible variety $X$ over $k$. Assume that $\dim(p_J\times p_J)(f(X))>\sum_{i\in J}n_i$ for every non-empty subset $J$ of $\{1,\ldots,s\}$ $($in particular, $\dim f(X)>\sum_{i=1}^sn_i)$. Then $f^{-1}(\Delta)$ is connected.\end{theorem}

Theorem \ref{debarre} actually asserts that the pair $(f(X),f(X)\cap\Delta)$ satisfies condition (1) of Theorem 
\ref{(2.1.24)} below. Then by Theorem \ref{(2.1.24)} the pair $(f(X),f(X)\cap\Delta)$ also satisfies condition (2) of Theorem \ref{(2.1.24)}, i.e. $K(f(X))$ is an algebraically closed subfield of $K(f(X)_{/f(X)\cap\Delta})$. 

The aim of this paper is to prove, under a slightly stronger hypothesis,  the following strengthening of Theorem \ref{debarre} of Debarre:

\begin{theorem}\label{main} Under the above notation, let $f\colon X\to P\times P$ be a morphism from a complete irreducible variety $X$ over $k$. Assume that $\dim(p_J\times p_J)(f(X))>\sum_{i\in J}n_i+p-1$ for every non-empty subset $J$ of $\;\{1,\ldots,s\}$ with $p$ elements $($and in particular, $\dim f(X)>\sum_{i=1}^sn_i+s-1)$. Then $f^{-1}(\Delta)$ is $G3$ in $X$.\end{theorem}

Theorem \ref{main} generalizes the main result of \cite{B} (see also \cite{B1}, chap. 10, which corresponds to the case $s=1$, i.e. when $P$ is a projective space). In Section 4 we give some relevant consequences of Theorem \ref{main} (see Corollaries \ref{debarre2},  \ref{debarre3}, \ref{debarre4} and \ref{last}). For instance one of them asserts that, if $Z$ is a closed irreducible subvariety of $\mathbb P^{n_1}\times\mathbb P^{n_2}$ of codimension 
$<\frac{1}{2}\min\{n_1,n_2\}$, then the diagonal $\Delta_Z$ of $Z\times Z$ is $G3$ in $Z\times Z$. 

We do not know if in Theorem \ref{main} the hypothesis that $\dim(p_J\times p_J)(f(X))>\sum_{i\in J}n_i+p-1$ can be relaxed to $\dim(p_J\times p_J)(f(X))>\sum_{i\in J}n_i$ for every non-empty subset $J$ of $\{1,\ldots,s\}$ with $p$ elements (as in Theorem \ref{debarre} of Debarre).

The proof of Theorem \ref{main} is given in Section 4 and makes use of the so-called join construction together with a systematic use of some basic results on formal-rational functions (which are recalled in Section 2). One of the main ingredients of the proof of Theorem \ref{main} is the following result (which is proved in Section 3):

\begin{theorem}\label{main1}  Under the above notation, let $f\colon X\to P= \mathbb P^{n_1}\times\cdots\times\mathbb P^{n_s}$ be a morphism from a complete irreducible variety $X$ over $k$. Assume that for every $i=1,\ldots,s$ we are given a linear subspace $L_i$ of $\;\mathbb P^{n_i}$ of codimension $r_i>0$ such that $\dim p_J(f(X))>\sum_{i\in J}r_i+p-1$ for every non-empty subset $J$ of $\;\{1,\ldots,s\}$ with $p$ elements $($and in particular, $\dim f(X)>\sum_{i=1}^sr_i+s-1)$. Then $f^{-1}(L_1\times\cdots\times L_s)$ is $G3$ in $X$.
\end{theorem}

As in Theorem \ref{main}, we do not know whether the hypothesis that $\dim p_J(f(X))>\sum_{i\in J}r_i+p-1$ can be relaxed to $\dim p_J(f(X))>\sum_{i\in J}r_i$ (as in Theorem \ref{debarre0}, (2) of Debarre below). We want to stress that the proofs Theorems \ref{main} and \ref{main1} make essential use of Debarre's connectivity results (Theorems 
\ref{debarre} above and \ref{debarre0} below) together with a number of results on formal functions recalled in Section 2. 

Theorem \ref{main1} (which, as any Bertini-type result, may also have an interest in its own) represents a generalization to a product of projective spaces of the following important result due to Grothendieck and Faltings:

\begin{theorem}[Grothendieck--Faltings \cite{F}]\label{gf} Let $X$ be a closed irreducible subvariety of $\mathbb P^n$ of dimension $d\geq 2$, and let $Y$ be the set-theoretic intersection of $X$ with $r$ hyperplanes of
$\mathbb P^n$, with $1\leq r\leq d-1$. Then $Y$ is $G3$ in $X$.\end{theorem}

Notice that in the above theorem Grothendieck had first proved the fact that $Y$ is connected (see \cite{SGA2}, Expos\'e XIII, Corollaire 2.3, for an even more general and a slightly stronger result), while later on Faltings improved Grothendieck's result to get the stronger conclusion that $Y$ is $G3$ in $X$ (see \cite{F}). Grothendieck and Faltings used local methods in their proofs. However, such local methods do not seem appropriate to prove Theorem \ref{main1} above. Therefore one has to appeal to global geometric methods. And this is done in Section 3 by showing that one of the global proofs of Theorem \ref{gf}  given by P. Bonacini, A. Del Padrone and M. Nesci in \cite{BDN}, suitably modified (and making also use of the results on formal functions recalled in the next section), works in our new situation as well. 

\medskip

Throughout this paper we shall fix an algebraically closed ground field $k$ of arbitrary characteristic.
The terminology and notation used are standard, unless otherwise specified.

\medskip

{\bf Acknowledgment.} We want to thank the referee for some suggestions which led to a slight improvement of the presentation.


\section{Background material}\label{first}\addtocounter{subsection}{1}

In this section we gather together the known results that are going to be used in Sections 3 and 4.

\begin{theorem}[Hironaka-Matsumura \cite{HM}, or also \cite{B1}, Thm. 9.11]\label{(2.1.12)} Let $f\colon X'\to X$ be
a proper surjective morphism of irreducible varieties over $k$. Then for every closed
subvariety $Y$ of $X$ there is a canonical isomorphism
$$K({X'}_{/f^{-1}(Y)})\cong[K(X')\otimes_{K(X)}K({X}_{/Y})]_0,$$
where $[A]_0$ denotes the
total ring of fractions of a commutative unitary ring $A$.\end{theorem}

\begin{corollary}\label{(2.1.14)}  Under the hypotheses of Theorem
$\ref{(2.1.12)}$, assume that $\,Y$ is $G3$ in $X$. Then $f^{-1}(Y)$ is $G3$ in
$X'$.\end{corollary}

\begin{prop}[Hironaka--Matsumura \cite{HM}, or also \cite{B1}, Cor. 9.10]\label{(2.1.11)} Let $X$ be an irreducible algebraic variety over $k$,
and let $Y$ be closed subvariety of $X$. Let $u\colon X'\to X$ be the $($birational$)$
normalization of $X$. Then $K(X_{/Y})$ is a field if and only if
$u^{-1}(Y)$ is connected. \end{prop}

\begin{prop}[\cite{B1}, Prop. 9.23]\label{(2.1.25)}  Let $f\colon X'\to X$ be a proper surjective
morphism of irreducible algebraic varieties over $k$, and let $Y\subset X$ and
$Y'\subset X'$ be closed subvarieties such that $f(Y')\subseteq Y$.
Assume that the rings $K(X_{/Y})$  and  $K({X'}_{/f^{-1}(Y)})$ are both fields.
If $Y'$ is $G3$ in $X'$, then $Y$ is $G3$ in $X$.\end{prop}

\begin{theorem}[Hironaka--Matsumura \cite{HM}]\label{(2.1.222)} Let $n_1,\ldots,n_s$ be $s\geq 2$ positive integers, and let $N_i$ be a line in the projective space $\mathbb P^{n_i}$ over $k$, $i=1,\ldots,s$. Then $N_1\times\ldots\times N_s$ is $G3$ in $\mathbb P^{n_1}\times\cdots\times\mathbb P^{n_s}$.\end{theorem}

This result follows easily from \cite{HM}, Lemma (3.1) and Theorem (2.13). Notice also that, if $s=2$, Theorem \ref{(2.1.222)} is a special case of a subsequent more general result of Speiser \cite{Sp} asserting that every irreducible subvariety $Y$ of $\mathbb 
P^{n_1}\times\mathbb P^{n_2}$, such that $\dim p_i(Y)>0$ for $i=1,2$, is $G3$ in $\mathbb P^{n_1}\times\mathbb P^{n_2}$, where $p_1$ and $p_2$ are the canonical projections of $\mathbb P^{n_1}\times\mathbb P^{n_2}$ (in fact, with some obvious modifications, Speiser's result can be generalized for a product of $s$ projective spaces). In Section 2 we are going to use  only the special case stated in Theorem \ref{(2.1.222)}.

\medskip

Let now $f\colon X'\to X$ be a proper surjective morphism of irreducible varieties, and let $Y\subset X$
and $Y'\subset X'$ be closed subvarieties such that $f(Y')\subseteq Y$. Then one can define a canonical
map of $k$-algebras $\w f^*\colon K(X_{/Y})\to K(X'_{/Y'})$ (pull back of formal-rational functions, see \cite{HM}, or also \cite{B1}, Cor. 9.8) rendering commutative the following diagram:

\begin{equation*}
\begin{CD}
K(X)@>f^*>> K(X')\\
@V\alpha_{X,Y} VV @ VV\alpha_{X',Y'}V\\
K(X_{/Y})@ >\w f^*>> K(X'_{/Y'})\\
\end{CD}
\end{equation*}

\begin{theorem}[B\u adescu--Schneider \cite{BSch}, or also \cite{B1}, Thm. 9.21]\label{(2.1.23)} 
Let $\zeta\in K(X_{/Y})$ be a formal-rational function of an irreducible
variety $X$ over $k$ along a closed subvariety $Y$ of $X$ such that
$K({X}_{/Y})$ is a field. Then the following conditions are equivalent:
\begin{enumerate}
\item[(1)] $\zeta$ is algebraic over $K(X)$.
\item[(2)] There is a proper surjective morphism $f\colon X'\to X$ from an
irreducible variety $X'$ and a closed
subvariety $Y'$ of $X'$ such that $f(Y')\subseteq Y$ and
$\w{f}^*(\zeta)$ is algebraic over $K(X')$.
\item[(3)] There is a proper surjective morphism $f\colon X'\to X$ from
an irreducible variety $X'$ and a closed
subvariety $Y'$ of $X'$ such that $f(Y')\subseteq Y$ and
$\tilde{f}^*(\zeta)\in K(X')$ $($more precisely, there exists a rational
function $t\in K(X')$ such that $\w f^*(\zeta)=\alpha_{X',Y'}(t))$.
\end{enumerate}\end{theorem}

\begin{theorem}[B\u adescu--Schneider \cite{BSch}, or also \cite{B1}, Cor. 9.22]\label{(2.1.24)}  Let $(X,Y)$ be a pair consisting of a complete irreducible variety $X$ over $k$ and
a closed subvariety $Y$ of $X$. The following conditions are
equivalent:
\begin{enumerate}
\item[(1)] For every proper surjective morphism $f\colon X'\to X$
from an irreducible variety $X'$, $f^{-1}(Y)$ is connected.
\item[(2)] $K({X}_{/Y})$ is a field and $K(X)$ is algebraically
closed in $K({X}_{/Y})$.
\end{enumerate}
\end{theorem}

Finally we shall make use of the following two theorems of Debarre:

\begin{theorem}[Debarre \cite{De}]\label{debarre0} Under the notation of the introduction, let $f\colon X\to P=\mathbb P^{n_1}\times\cdots\times\mathbb P^{n_s}$ be a morphism from a complete irreducible variety $X$ over $k$, and assume that for every $i=1,\ldots,s$ we are given a linear subspace $L_i$ of $\;\mathbb P^{n_i}$ of codimension $\geq 1$.
\begin{enumerate}
\item[(1)] If $\dim p_J(f(X))\geq\sum_{i\in J}\codim_{\mathbb P^{n_i}}(L_i)$ for every non-empty subset $J$ of $\{1,\ldots,s\}$, then $f^{-1}(L_1\times\cdots\times L_s)\neq\varnothing$.
\item[(2)] If $\dim p_J(f(X))>\sum_{i\in J}\codim_{\mathbb P^{n_i}}(L_i)$ for every non-empty subset $J$ of $\{1,\ldots,s\}$, then $f^{-1}(L_1\times\cdots\times L_s)$ is connected.
\end{enumerate}\end{theorem}

\section{Proof of Theorem \ref{main1}}

In this section we prove Theorem \ref{main1} of the introduction. For the simplicity of notation we shall give the proof in the case of a product of two projective spaces (the proof in general goes similarly, with only small modifications).
Fix a product of projective spaces $\mathbb P^{n_1}\times\mathbb P^{n_2}$ over $k$, with $n_i\geq 1$, $i=1,2$, and denote by $p_1$ and $p_2$ the canonical projections of $\mathbb P^{n_1}\times\mathbb P^{n_2}$. Then Theorem \ref{main1} becomes:

\begin{theorem}\label{faltings1}  Under the above notation, let $X$ be a projective irreducible variety, and let $f\colon X\to\mathbb P^{n_1}\times\mathbb P^{n_2}$ be a morphism. Assume that for $i=1,2$ we are given a linear subspace $L_i$ of $\;\mathbb P^{n_i}$ of codimension $r_i>0$ such that $\dim f(X)>r_1+r_2+1$, $\dim p_1(f(X))>r_1$ and $\dim p_2(f(X))>r_2$. Then $f^{-1}(L_1\times L_2)$ is $G3$ in $X$.\end{theorem}

\proof To prove the theorem we show that one of the global proofs of Theorem \ref{gf} of Grothendieck--Faltings given by P. Bonacini, A. Del Padrone and M. Nesci in \cite{BDN}, suitably modified (and making free use of the results on formal functions recalled in the previous section), works in our new situation as well. 

First of all, by Corollary \ref{(2.1.14)} above, we can replace $X$ by $f(X)$ and $f^{-1}(L)$ by $f(X)\cap L$, where $L:=L_1\times L_2$. In  other words, we may assume that $X$ is a closed irreducible subvariety of $P:=\mathbb P^{n_1}\times\mathbb P^{n_2}$ such that $\dim X> r_1+r_2+1$, $\dim p_1(X)>r_1$ and $\dim p_2(X)>r_2$, and then we have to prove that $Y:=X\cap L$ is $G3$ in $X$. 

Let $s_1=\cdots=s_{r_1}=0$ (resp. $t_1=\cdots=t_{r_2}=0$) be equations defining $L_1$ in $\mathbb P^{n_1}$ (resp.  $L_2$ in $\mathbb P^{n_2}$), where  $s_i\in H^0(\mathscr O_{\mathbb P^{n_1}}(1))$ and $t_j\in H^0(\mathscr O_{\mathbb P^{n_2}}(1))$. Let
$$\alpha_1\colon H^0(\mathscr O_{\mathbb P^{n_1}}(1))\to H^0(\mathscr O_{\mathbb P^{n_1}\times\mathbb P^{n_2}}(1,0))\;\text{and}$$
$$\alpha_2\colon H^0(\mathscr O_{\mathbb P^{n_2}}(1))\to H^0(\mathscr O_{\mathbb P^{n_1}\times\mathbb P^{n_2}}(0,1))$$
be the canonical injective maps, and set $V_i:=\Im(\alpha_i)$, $i=1,2$. In particular, if $\sigma=\alpha_1(\sigma')\in V_1$ and $x=(x_1,x_2)\in \mathbb P^{n_1}\times\mathbb P^{n_2}$, one has $\sigma(x)=\sigma'(x_1)$. If we set 
$$Q:=\mathbb P(V_1^{\oplus(r_1+1)})\times\mathbb P(V_2^{\oplus(r_2+1)}),$$
consider the incidence variety 
$$Z:=\{(x,[\sigma_0,\ldots,\sigma_{r_1}],[\tau_0,\ldots,\tau_{r_2}])\in X\times Q\,|\,\sigma_i(x)=
\tau_j(x)=0,$$
$$\forall i=0,\ldots,r_1,\,\forall j=0,\,\ldots r_2\}.$$
Clearly $Z$ is a closed subset of $X\times Q$ and, in particular, $Z$ is a projective variety. Denote by 
$g\colon Z\to X$ and $h\colon Z\to Q$ the restrictions to $Z$ of the canonical projections of $X\times Q$. Then $g$ and $h$ are proper morphisms. Let $([\sigma_0,\ldots,\sigma_{r_1}],[\tau_0,\ldots,\tau_{r_2}])\in Q$ be an arbitrary point, where $\sigma_i=\alpha_1(\sigma'_i)$ and $\tau_j=\alpha_2(\tau'_j)$, $i=0,\ldots,r_1$, $j=0,\ldots,r_2$ ($\sigma'_i\in H^0(\mathscr O_{\mathbb P^{n_1}}(1))$ and
$\tau'_j\in H^0(\mathscr O_{\mathbb P^{n_2}}(1))$).  Then 
$$h^{-1}([\sigma_0,\ldots,\sigma_{r_1}],[\tau_0,\ldots,\tau_{r_2}])\cong X\cap(M_1\times M_2), $$
with $M_1=\{y\in\mathbb P^{n_1} | \sigma'_i(y)=0, i=0,\ldots,r_1\}$ and $M_2=\{z\in\mathbb P^{n_2}|\tau'_j(z)=0, j=0,\ldots,r_2\}$. Since $\dim p_1(X)\geq r_1+1$, $\dim p_2(X)\geq r_2+1$ and  $\dim X\geq r_1+r_2+2=(r_1+1)+(r_2+1)$, by Theorem \ref{debarre0}, (1)  we infer that $X\cap(M_1\times M_2)\neq\varnothing$. It follows  that the morphism $h$ is surjective.

On the other hand, the fibers of $g$ are all isomorphic to the product $\mathbb P^{n_1(r_1+1)-1}\times\mathbb P^{n_2(r_2+1)-1}$. In particular, the morphism $g$ is also surjective. Moreover, since  $X$ is projective and irreducible and all the fibers of $g$ are projective, irreducible and of the same dimension, by an elementary result  (see   \cite{Sh}, Part 1, p. 77, Theorem 8) it follows that $Z$ is also irreducible.

Notice that for every $([\sigma_0,\ldots,\sigma_{r_1}],[\tau_0,\ldots,\tau_{r_2}])\in Q$ such that 
$\sigma'_0,\ldots,\sigma'_{r_1}\in H^0(\mathscr O_{\mathbb P^{n_1}}(1))$ generate the same subspace as $s_1,\ldots,s_{r_1}$ and
$\tau'_0,\ldots,\tau'_{r_2}\in H^0(\mathscr O_{\mathbb P^{n_2}}(1))$ generate the same subspace as $t_1,\ldots,t_{r_2}$ (where $\sigma_i=\alpha_1(\sigma'_i)$, $i=0,\ldots,r_1$, and
$\tau_j=\alpha_2(\tau'_j)$, $j=0,\ldots,r_2$), we have
\begin{equation}\label{e:xx}h^{-1}([\sigma_0,\ldots,\sigma_{r_1}],[\tau_0,\ldots,\tau_{r_2}])\subset g^{-1}(Y).\end{equation}
Let $N_1$ denote the line of $\mathbb P(H^0(\mathscr O_{\mathbb P^{n_1}}(1)^{\oplus(r_1+1)})\cong\mathbb P(V_1^{\oplus(r_1+1)})$ passing through the points $[s_1,\ldots,s_{r_1},0]$ and $[0,s_1,\ldots,s_{r_1}]$, and $N_2$  the line of $\mathbb P(H^0(\mathscr O_{\mathbb P^{n_2}}(1)^{\oplus(r_2+1)})\cong\mathbb P(V_2^{\oplus(r_2+1)})$ passing through the points $[t_1,\ldots,t_{r_2},0]$ and $[0,t_1,\ldots,t_{r_2}]$.  By Theorem \ref{(2.1.222)} the product 
$N_1\times N_2$ is $G3$ in $Q=\mathbb P(V_1^{\oplus(r_1+1)})\times\mathbb P(V_2^{\oplus(r_2+1)})$, whence by Corollary 
\ref{(2.1.14)}, $h^{-1}(N_1\times N_2)$ is also $G3$ in $Z$ (because $h$ is proper and surjective). On the other hand, since every point of the line $N_1$ is of the form $[\lambda s_1,\lambda s_2+\mu s_1,\ldots,\lambda s_{r_1}+\mu s_{r_1-1},\mu s_{r_1}]$, with $(\lambda,\mu)\in k^2\setminus\{(0,0)\}$ (and similarly for the line $N_2$), by $\eqref{e:xx}$ we get the inclusion
$$h^{-1}(N_1\times N_2)\subset g^{-1}(Y),$$
and therefore the natural maps of $k$-algebras
$$K(Z)\to K(Z_{/g^{-1}(Y)})\to K(Z_{/h^{-1}(N_1\times N_2)}),$$
whose composition is an isomorphism because $h^{-1}(N_1\times N_2)$ is $G3$ in $Z$. We claim that this implies that the first map is also an isomorphism. To see this, it will be enough to check that $K(Z_{/g^{-1}(Y)})$ is a field. This latter fact follows from Debarre's Theorem \ref{debarre0}, (2) and from Proposition \ref{(2.1.11)}.
Indeed, if $u\colon \tilde Z\to Z$ is the birational normalization of $Z$, Theorem \ref{debarre0}, (2) together with the surjectivity of $g$ imply that $u^{-1}(g^{-1}(Y))$ is connected. Then by Proposition 
\ref{(2.1.11)}, the ring 
$K(Z_{/g^{-1}(Y)})$ of formal-rational functions of $Z$ along $g^{-1}(Y)$ is a field. 

Therefore we proved that $g^{-1}(Y)$ is $G3$ in $Z$. Moreover, as above using Theorem 
\ref{debarre0}, (2) we see that $K(X_{/Y})$ is a field. Since $g$ is proper and surjective, by Proposition \ref{(2.1.25)}  we finally get the fact that $Y$ is $G3$ in $X$.\qed

\begin{rem}\label{f2} {\em The conclusion of Theorem \ref{main1} (and in particular, of Theorem \ref{faltings1}) remains still true if one only assume that there exists a Zariski open subset $U$ of $\mathbb P^{n_1}\times\cdots\times\mathbb P^{n_s}$ such that $L\subset U$, $f(X)\subset U$, and the corestriction $f\colon X\to U$ is a proper morphism. Indeed, according to the proof of Theorem \ref{main1} one may assume that $X\subset U$ and $X$ is closed in $U$. Let $\overline X$ be the closure of $X$ in $P$. Since $L\subset U$, $X\cap L=\overline X\cap L$. Then by Theorem \ref{main1} applied to $\overline X$, $Y=\overline X\cap L$ is $G3$ in $\overline X$, i.e. the natural map $K(\overline X)\to K(\overline X_{/Y})$ is an isomorphism. But clearly $K(\overline X)=K(X)$  and $K(\overline X_{/Y})=K(X_{/Y})$ (because $X$ is open in $\overline X$).}\end{rem}

\section{Proof of Theorem \ref{main} and some consequences}

Using Theorem \ref{main1} (resp. Theorem \ref{faltings1} of the previous section), we can prove Theorem \ref{main}. As in the previous section we shall give the proof in the case $s=2$, in other words, we shall prove the following:

\begin{theorem}\label{debarre1} Let $f\colon X\to P\times P$ be a morphism from a complete irreducible variety $X$ over $k$, where $P:=\mathbb P^{n_1}\times\mathbb P^{n_2}$, with $n_1,n_2\geq 1$, such that $\dim f(X)>n_1+n_2+1$, $\dim (p_1\times p_1)(f(X))>n_1$ and $\dim (p_2\times p_2)(f(X))>n_2$ $($with $p_1$ and $p_2$ the canonical projections of $\mathbb P^{n_1}\times\mathbb P^{n_2})$. Then $f^{-1}(\Delta)$ is $G3$ in $X$.
\end{theorem}

\proof As in the proof of Theorem \ref{faltings1}, we may assume that $X\subseteq P\times P$ and $f$ the canonical inclusion.  Set:
$$U_{n_1}:=\mathbb P^{2n_1+1}\setminus (L_1^{n_1}\cup L_2^{n_1})\:\:\text{and}\;\;U_{n_2}:=\mathbb 
P^{2n_2+1}\setminus (L_1^{n_2}\cup L_2^{n_2}),$$
where $L_1^{n_1}$ (resp. $L_2^{n_1}$) is the linear subspace of $\mathbb P^{2n_1+1}$ (of homogeneous coordinates $[t_0,\ldots,t_{n_1},t'_0,\ldots,t'_{n_1}]$)
of equations $t_0=\cdots=t_{n_1}=0$ (resp. $t'_0=\cdots=t'_{n_1}=0$), and $L_1^{n_2}$ (resp. $L_2^{n_2}$) is the linear subspace of $\mathbb P^{2n_2+1}$ (of homogeneous coordinates $[u_0,\ldots,u_{n_2},u'_0,\ldots,u'_{n_2}]$) of equations $u_0=\cdots =u_{n_2}=0$ (resp. $u'_0=\cdots =u'_{n_2}=0$).
Since $t_i-t_i'$ and $u_j-u_j'$ are homogeneous elements  of degree $1$, $i=0,\ldots ,n_1$, $j=0,\ldots ,n_2$, it makes also sense to consider the linear subspace
$H_{n_1}$ of $\mathbb P^{2n_1+1}$ of equations 
$$t_0-t_0'=\cdots=t_{n_1}-t_{n_1}'=0,$$ 
and the linear subspace  $H_{n_2}$ of
$\mathbb P^{2n_2+1}$ of equations 
$$u_0-u_0'=\cdots=u_{n_2}-u_{n_2}'=0.$$
Clearly $H_{n_1}\subset U_{n_1}$ and $H_{n_2}\subset U_{n_2}$, whence  
$$H_{n_1}\times H_{n_2}\subset U:=U_{n_1}\times U_{n_2}\subset \mathbb P^{2n_1+1}\times \mathbb P^{2n_2+1}.$$
Consider the rational map
$$\pi\colon\mathbb P^{2n_1+1}\times \mathbb P^{2n_2+1}\dashrightarrow P\times P=\mathbb P^{n_1}\times \mathbb P^{n_2}\times\mathbb P^{n_1}\times \mathbb P^{n_2}$$
defined by  
$$\pi([t_0,\ldots ,t_{n_1}, t_0',\ldots ,t_{n_1}'],[u_0,\ldots ,u_{n_2},u_0'\ldots ,u_{n_2}'])=$$
$$= ([t_0,\ldots ,t_{n_1}],[u_0,\ldots ,u_{n_2}],[ t_0',\ldots ,t_{n_1}'],[u_0',\ldots ,u_{n_2}']),$$
$\forall\ ([t_0,\ldots ,t_{n_1}, t_0',\ldots ,t_{n_1}'],[u_0,\ldots ,u_{n_2},u'_0\ldots ,u_{n_2}'])\in\mathbb P^{2n_1+1}\times \mathbb P^{2n_2+1}$.
Actually, modulo the canonical isomorphism 
$$\mathbb P^{n_1}\times \mathbb P^{n_2}\times\mathbb P^{n_1}\times \mathbb P^{n_2}\cong\mathbb P^{n_1}\times \mathbb P^{n_1}\times\mathbb P^{n_2}\times \mathbb P^{n_2},$$
$\pi$ is nothing but the product $\pi_1\times\pi_2$ of the rational maps $\pi_1\colon\mathbb P^{2n_1+1}\dasharrow \mathbb P^{n_1}\times \mathbb P^{n_1}$ and 
$\pi_2\colon\mathbb P^{2n_2+1}\dasharrow \mathbb P^{n_2}\times \mathbb P^{n_2}$ defined by 
$$\pi_1([t_0,\ldots,t_{n_1},t'_0,\ldots,t'_{n_1}])=([t_0,\ldots,t_{n_1}],[t'_0,\ldots,t'_{n_1}]),$$
$$\pi_2([u_0,\ldots,u_{n_2},u'_0,\ldots,u'_{n_2}])=([u_0,\ldots,u_{n_2}],[u'_0,\ldots,u'_{n_2}]).$$
Then  the  map $\pi$  is defined precisely on the open subset $U=U_{n_1}\times U_{n_2}$ and $\pi$ is the projection of a locally trivial $(\mathbb G_m\times\mathbb G_m)$-bundle, where $\mathbb G_m=k\setminus\{0\}$; in particular, all the fibers of $\pi$ are isomorphic to $\mathbb{G}_m\times\mathbb{G}_m$. It is clear that the restriction map 
$\pi |(H_{n_1}\times H_{n_2})$ defines an isomorphism $H_{n_1}\times H_{n_2}\cong \Delta$.

Now consider the commutative diagram
\medskip
\begin{equation*}
\begin{CD}
Y':=U_{X}\cap(H_{n_1}\times H_{n_2})@>>> U_{X}@>>>U\\
@V\cong VV @ V\pi_{X}VV@ VV\pi V\\
Y:=X\cap \Delta @ >>> X@>f>>P\times P\\
\end{CD}
\end{equation*}

\medskip

\noindent where $U_{X}:=\pi^{-1}(X)$ and $\pi_{X}:=\pi|U_{X}$ (the restriction of $\pi$ to $U_{X}$), and the horizontal arrows are closed embeddings.

Set $X_{13}:=(p_1\times p_1)(X)$ and $X_{24}:=(p_2\times p_2)(X)$, and let $q_1$ and $q_2$ be the canonical projections of $\mathbb P^{2n_1+1}\times\mathbb P^{2n_2+1}$. Since $\pi=\pi_1\times\pi_2$ and $U_X=\pi^{-1}(X)$ it follows that $q_1(U_X)=U_{X_{13}}$ and $q_2(U_X)=U_{X_{24}}$, where
$U_{X_{13}}:=\pi_1^{-1}(X_{13})$ and $U_{X_{24}}:=\pi_2^{-1}(X_{24})$ ($U_{X_{13}}$ and $U_{X_{24}}$ are locally trivial $\mathbb G_m$-bundles over $X_{13}$ and $X_{24}$ respectively). It follows that 
$$\dim q_1(U_X)=\dim U_{X_{13}}=\dim X_{13}+1>n_1+1, \;\text{and}$$
$$\dim q_2(U_X)=\dim U_{X_{24}}=\dim X_{24}+1>n_2+1,$$
where the last two inequalities follow from the hypotheses. Moreover, by hypothesis we also have
$$\dim U_X=\dim X+2>n_1+n_2+3=(n_1+1)+(n_2+1)+1.$$
Thus, if we take $r_1=n_1+1$, $r_2=n_2+1$ we may apply Theorem \ref{faltings1} (via Remark \ref{f2})
to deduce that $Y'$ is $G3$ in $U_X$. 

Let $W$ be the closure of $U_X$ in $\mathbb P^{2n_1+1}\times\mathbb P^{2n_2+1}$, and let $W'$ be the graph of the rational map $\pi_X\colon W\dasharrow X$ (i.e. the closure in $W\times X$ of the graph of $\pi_X\colon U_X\to X$). Then we get the morphisms
$u\colon W'\to W$ and $v\colon W'\to X$, with $v$ surjective, such that $\pi_X\circ u=v$. Moreover, $W$ contains $U_X$ as an open subset such that $v|U_X=\pi_X$ and $v(Y')=Y$. Then from Theorem \ref{(2.1.23)} it follows that the field extension given by $\alpha_{X,Y}\colon K(X)\to K(X_{/Y})$ is algebraic.
On the other hand, Debarre's Theorem \ref{debarre} asserts that the pair $(X,Y)$ satisfies condition (1) of Theorem 
\ref{(2.1.24)}, whence by this latter theorem $K(X)$ is also algebraically closed in $K(X_{/Y})$ (via the map $\alpha_{X,Y}\colon K(X)\to K(X_{/Y})$).
It follows that the map $\alpha_{X,Y}$ is an isomorphism, i.e. $Y$ is $G3$ in $X$.\qed

\begin{corollary}\label{debarre2} Let $X$ be a closed irreducible variety  of $P\times P$, with $P=\mathbb P^{n_1}\times\mathbb P^{n_2}$, and $n_1\geq n_2\geq 1$. Assume that $\codim_{P\times P}X<n_2$. Then $X\cap\Delta$ is $G3$ in $X$.\end{corollary}

\proof We show that the hypotheses of Theorem \ref{debarre1} are fulfilled for the embedding of $X$ in $P\times P$. Clearly, $\dim X\geq 2n_1+n_2+1>n_1+n_2+1$. We also have to prove that $\dim X_{13}>n_1$ and $\dim X_{24}>n_2$. Assuming for example that
$\dim X_{13}\leq n_1$, we see that the general fiber $F$ of  $X\to X_{13}$ has dimension 
$\dim X-\dim X_{13}\geq 2n_1+n_2+1-\dim X_{13}\geq (2n_1+n_2+1)-n_1=n_1+n_2+1\geq 2n_2+1$. But this is absurd because $F\cong(p_2\times p_2)(F)\subseteq \mathbb P^{n_2}\times\mathbb P^{n_2}$. The other inequality follows similarly. Then by Theorem \ref{debarre1} we get the conclusion.\qed

\begin{corollary}\label{debarre3} Let $Y$ and $Z$ be two closed irreducible subvarieties of $P=\mathbb P^{n_1}\times\mathbb P^{n_2}$, with $n_1\geq n_2\geq 1$. Assume that $\dim Y+\dim Z>2n_1+n_2$.
Then $Y\cap Z\cong (Y\times Z)\cap\Delta$ is $G3$ in $Y\times Z$, where $\Delta$ is the diagonal of $P\times P$. In particular, if $\codim_PZ<\frac{n_2}{2}$ then the 
diagonal $\Delta_Z$ of $Z\times Z$ is $G3$ in $Z\times Z$.\end{corollary}

\proof Apply Corollary \ref{debarre2} to $X=Y\times Z$.\qed

\begin{corollary}\label{debarre4} Let $Y$ and $Z$ be two closed irreducible subvarieties of $P=\mathbb P^{n_1}\times\mathbb P^{n_2}$, with $n_1\geq n_2\geq 1$. Assume that $\dim Y+\dim Z>2n_1+n_2$.
Then $Y\cap Z$ is $G3$ in $Y$ and in $Z$. \end{corollary}

\proof If $p\colon Y\times Z\to Y$ is the first projection, we clearly have $p((Y\times Z)\cap\Delta)=Y\cap Z$. Since $\dim Y+\dim Z>2n_1+n_2$, we can apply Corollary \ref{debarre3} to deduce that $(Y\times Z)\cap\Delta\cong Y\cap Z$ is $G3$ in $Y\times Z$. If  $f\colon V\to Y$ is an arbitrary proper surjective morphism. By Corollary \ref{(2.1.14)} it follows that $(f\times\id_Z)^{-1}(\Delta)$ is $G3$ in $V\times Z$, and in particular, $(f\times\id_Z)^{-1}(\Delta)\cong f^{-1}(Y\cap Z)$ is connected.  Then by Proposition \ref{(2.1.11)} above it follows that $K(Y_{/Y\cap Z})$ is a field. Moreover, the same argument shows that for every proper surjective morphism $g\colon U\to Y\times Z$ (from an irreducible variety 
$U$), $g^{-1}(p^{-1}(Y\cap Z))$ is connected, whence  
$K((Y\times Z)_{/p^{-1}(Y\cap Z)})$ is also a field (by Proposition \ref{(2.1.11)} again).

Since  $(Y\times Z)\cap\Delta$ is $G3$ in $Y\times Z$,  by Proposition \ref{(2.1.25)} applied to $p\colon Y\times Z\to Y$, we infer that  $Y\cap Z$ is $G3$ in $Y$. Similarly one proves that $Y\cap Z$ is $G3$ in $Z$. \qed

\begin{rem}\label{s233} {\em   Corollary \ref{debarre4} is an analogue of the following result of Faltings \cite{F}: for every closed irreducible subvarieties $Y$ and $Z$ of $P=\mathbb P^n$ such that 
$\dim Y+\dim Z>n$, $Y\cap Z$ is $G3$ in $Y$ and in $Z$. }\end{rem}

To state the last corollary we need to recall the following:

\begin{definition}\label{(2.2.1)} {\em Let $Y$ be a closed subvariety 
of a projective variety $X$. According to Grothendieck \cite{SGA2}, Expos\'e X, one says  that
the pair $(X,Y)$ satisfies the  Grothendieck--Lefschetz condition
$\Lef(X,Y)$ if for every open subset $U$ of $X$ containing $Y$ and for every vector
bundle $E$ on $U$ the natural map $H^0(U,E)\to H^0(X_{/Y},\h E)$ is an
isomorphism, where $\h E=\pi^*(E)$, with $\pi\colon X_{/Y}\to U$ the canonical
morphism. We also say that $(X,Y)$ satisfies the {\it effective
Grothendieck-Lefschetz condition} $\Leff(X,Y)$ if  $\Lef(X,Y)$ holds and, moreover, for every formal vector bundle $\mathcal E$ on $X_{/Y}$ there exists an open subset $U$ of $X$ containing $Y$ and a vector
bundle $E$ on $U$ such that $\mathcal E\cong\h E$.}\end{definition}

Then exactly as in the case of small--codimensional subvarieties of $\mathbb P^n$ (see \cite{B}, or also \cite{B1}, Theorem 11.7, p. 128), using Corollary \ref{debarre3}, one can prove the following:

\begin{corollary}\label{last} Let $Z$ be a closed irreducible subvariety of $\mathbb P^{n_1}\times\mathbb P^{n_2}$, with $n_1\geq n_2\geq 1$. If $\codim_PZ<\frac{n_2}{2}$, then $\Lef(Z\times Z,\Delta_Z)$ is satisfied, where $\Delta_Z$ is the diagonal of $Z\times Z$.\end{corollary}

\begin{rems}\label{last1} {\em 
\begin{enumerate}
\item[(1)] Corollary \ref{last} is an analogue of the following result proved in \cite{B} (see also \cite{B1}, Thm. 11.7, p. 128): if $Z$ is a closed irreducible subvariety in $\mathbb P^n$ of dimension $>\frac{n}{2}$ then the Grothendieck--Lefschetz condition $\Lef(Z\times Z,\Delta_Z)$ is satisfied. Notice that, in the case when the characteristic of $k$ is positive and $Z$ is locally Cohen--Macaulay, this latter statement is an old result of Speiser \cite{Sp2}  (proved by completely different methods). 
\item[(2)] On the other hand, as in the case of submanifolds of $\mathbb P^n$ (see \cite{B1}, Prop. 11.9, p. 129), it is not difficult to show that, under the hypotheses of Corollary \ref{last}, the effective Grothendieck--Lefschetz condition $\Leff(Z\times Z,\Delta_Z)$ is never satisfied.\end{enumerate}}
\end{rems}

In closing this paper we want to rise an open question. To do this we first need the following:

\begin{definition}\label{ca}  {\em   Let $P$ be a projective rational homogeneous space. Then, according to Goldstein \cite{Go} one defines the coampleness of $P$ as follows. 
Since $P$ is a homogeneous space, the tangent bundle $T_P$ of $P$ is generated by its global
sections; this implies that the tautological line bundle ${\mathcal O}_{{\mathbb P}(T_P)}(1)$ is also generated by its global sections. Then one defines the ampleness, $\amp(P)$, of $P$ as the maximum fiber dimension of the morphism 
$\varphi\colon{\mathbb P}(T_P)\to {\mathbb P}^N$ associated to the complete linear system $|{\mathcal
O}_{{\mathbb P}(T_P)}(1)|$. Finally, the coampleness, $\ca(P)$, of $P$ is defined  by $\ca(P):=\dim P-\amp(P)$. A result of Goldstein (\cite{Go}) asserts that $\ca(P)\geq r$, where $r$ is the minimum of ranks of the simple factors of the linear algebraic group $G$ acting transitively on $P$; in particular, $\ca(P)\geq 1$. For example, it is easy to see that $\ca(\mathbb P^n)=n$ (or more generally, $\ca( \mathbb{G}(d,n))=n$, where $\mathbb{G}(d,n)$ is the Grassmann variety of $d$-planes in $\mathbb P^n$) and $\ca(\mathbb P^{n_1}\times\mathbb P^{n_2})=n_2$, if $n_1\geq n_2\geq 1$.   }
\end{definition}

{\bf Question.}   Let $Z$ be a closed irreducible subvariety of a projective rational homogeneous space $P$ over $\mathbb C$. Is it true that if $\codim_PZ<\frac{1}{2}\ca(X)$, then the diagonal   $\Delta_Z$ of $Z\times Z$  is $G3$ in $Z\times Z$?

\medskip\medskip

The answer to this question is positive in the following cases:
\begin{enumerate}
\item If $Z=P$  (see \cite{BSch}, Theorem (4.16)).
\item If $P=\mathbb P^n$ (even in arbitrary characteristic, see \cite{B}, Corollary (3.1)).
\item If $P=\mathbb P^{n_1}\times\mathbb P^{n_2}$ (even in arbitrary characteristic, by Corollary  \ref{debarre3} above).
\end{enumerate}

\end{document}